\documentstyle[12pt]{article}
\begin{document}

\rightline{FTUV/95-45, IFIC/95-47}
\rightline{q-alg/9510011}
\rightline{To appear in J. Phys. A (1996)}

\begin{center}
{\Large {\bf Deformed  Minkowski
spaces: classification and properties}}
\end{center}

\begin{center}
{\large{J.A. de Azc\'{a}rraga$^{\dagger}$
and F. Rodenas$^{\dagger \; \ddagger}$}}
\end{center}

\noindent
{\small{\it $ \dagger$ Departamento de F\'{\i}sica Te\'{o}rica and IFIC,
Centro Mixto Universidad de Valencia-CSIC
E-46100-Burjassot (Valencia), Spain.}}

\noindent
{\small{\it $ \ddagger $ Departamento de Matem\'atica Aplicada,
 Universidad Polit\'ecnica de Valencia
E-46071 Valencia, Spain.}}

\vspace{1\baselineskip}

\begin{abstract}
Using general but simple covariance arguments, we classify
the `quantum'  Minkowski spaces for dimensionless deformation parameters.
This requires a previous analysis of the associated  Lorentz groups,
which reproduces    a previous classification
by  Woronowicz and Zakrzewski. As a consequence of the
unified analysis presented,  we  give the commutation
 properties, the deformed (and central) length element and  the
metric tensor for  the different  spacetime  algebras.
\end{abstract}

\section{Introduction}

\indent

Following the approach of  \cite{FTUV94-21},
we  present  here  a classification of the
possible deformed   Minkowski spaces (algebras).
Our analysis, which provides a common framework for
the properties of the various
Minkowski spacetimes,  requires the consideration of the two ($SL_q(2)$ and
 $SL_h(2)$)
deformations of
$SL(2,C)$  and provides a characterization  of the appropriate $R$-matrices
defining the  deformed Lorentz groups given in
\cite{WOZA} (see also \cite{PODWOR}).

It is well known that $GL(2,C)$  admits only
two different deformations having a central determinant:
one is the standard $q$-deformation \cite{DRINF,FRT1}
and the other is the non-standard or `Jordanian' $h$-deformation
\cite{DEMI,EOW,KAR}. Both  $GL_q(2)$ and
 $GL_h(2)$ have associated
`quantum spaces' in the sense of \cite{MANIN}. These
deformations (which may be shown to be related by contraction \cite{AKS})
are defined as  the associative algebras generated by the
entries  $a,b,c,d$ of a matrix $M$, the commutation properties of which
may be  expressed by an `FRT'  \cite{FRT1} equation
\begin{equation}\label{1.1}
R_{12}M_1M_2=M_2M_1R_{12}
\end{equation}
\noindent
for a suitable $R$-matrix. Let us summarize their properties.

\vspace{1\baselineskip}
\noindent
{\bf a)} For $GL_q(2)$ the $R$-matrix in (\ref{1.1})
is  ($\lambda \equiv q-q^{-1}$)
\begin{equation}\label{1.2}
R_q= \left[
\begin{array}{cccc}
q & 0 & 0 & 0 \\
0 & 1 & 0 & 0 \\
0 & \lambda & 1 & 0 \\
0 & 0 & 0 & q
\end{array} \right] \;, \;
\hat{R}_q \equiv {\cal P}R_q = \left[
\begin{array}{cccc}
q & 0 & 0 & 0 \\
0 & \lambda & 1 & 0 \\
0 & 1 & 0 & 0 \\
0 & 0 & 0 & q
\end{array} \right] \;, \; {\cal P}R_q {\cal P}=R^t_q \;,
\end{equation}
\noindent
where ${\cal P}$ is the permutation operator
(${\cal P}={\cal P}^{\dagger}$, ${\cal P}_{ij,kl}= \delta_{il} \delta_{jk}$),
and the commutation relations defining the quantum group algebra are
\begin{equation}\label{1.3}
\begin{array}{lll}
ab=qba \quad, & \quad  ac=qca \quad,  & \quad ad-da= \lambda bc \quad,\\
bc=cb \quad, & \quad bd=qdb \quad, & \quad cd=qdc \quad.
\end{array}
\end{equation}
\noindent
$Fun(GL_q(2))$ has   a quadratic central element,
\begin{equation}\label{1.4}
det_qM:=ad-qbc \quad;
\end{equation}
\noindent
$det_qM=1$  defines $SL_q(2)$.
The matrix $\hat{R}_q \equiv {\cal P} R_q$ satisfies Hecke's condition
\begin{equation}\label{1.4.1}
\hat{R}_q^{2} - \lambda \hat{R}_q - I = 0 \quad , \quad (\hat{R}_q - qI)
(\hat{R}_q + q^{-1}I) =0 \quad,
\end{equation}

\noindent
and (we shall assume  $q^2 \neq -1$ \cite{FRT1} throughout) it
has a spectral
decomposition in terms of  a rank three projector $P_{q+}$
and a rank one projector $P_{q-}$,
\begin{equation}\label{1.4.2}
\hat{R}_q = q P_{q+} -q^{-1} P_{q-} \;, \; \hat{R}_q^{-1} = q^{-1} P_{q+}
- q P_{q-}\;,\;
[\hat{R}_q , P_{q \pm} ] =0\; , \; P_{q \pm} \hat{R}_q P_{q \mp} = 0 \;,
\end{equation}
\begin{equation}\label{1.4.3}
P_{q +} =  \frac{I+q \hat{R}_q}{1+q^2} \quad , \quad
P_{ q -}=  \frac{I-q^{-1} \hat{R}_q}{1+q^{-2}}=
\frac{1}{1+q^{-2}}  \left[
\begin{array}{cccc}
0 & 0 & 0 & 0 \\
0 & q^{-2} & -q^{-1}  & 0 \\
0 & -q^{-1} & 1 & 0 \\
0 & 0 & 0 & 0
\end{array} \right]   \;.
\end{equation}
\noindent
The following relations  have an obvious  equivalent in the undeformed  case:
\begin{equation}\label{1.4.4}
\epsilon_qM^t\epsilon_q^{-1}=M^{-1} \;, \quad \epsilon_q= \left(
\begin{array}{cc}
      0 & q^{-1/2} \\
      -q^{1/2}&0
      \end{array} \right) = -\epsilon_q^{-1} \;, \quad
P_{q - \;ij,kl}=\frac{-1}{[2]_q}\epsilon_{q\;ij}\epsilon^{-1}_{q\;kl} \;.
\end{equation}
\noindent
The determinant of an ordinary $2 \times 2$ matrix  may be defined
as the proportionality
coefficient in  $(det M)P_{-} := P_{-} M_{1} M_{2}$ where $P_{-}$ is given by
(\ref{1.4.3}) for $q$=1. In  the  $q \neq 1$ case the $q$-determinant
(\ref{1.4}) may be expressed as
\begin{equation}\label{adet}
(det_q M)P_{q-} := P_{q-} M_{1} M_{2}\quad, \quad
(det_q M^{-1})P_{q-} = M_{2}^{-1} M_{1}^{-1}P_{q-} \quad,
\end{equation}
\noindent
($det_q M^{-1}= (det_q M)^{-1}$, and $(det_q M)^{\dagger}P_{q-}^{\dagger}
= M_{2}^{\dagger} M_{1}^{\dagger}P_{q-}^{\dagger}$).

\vspace{1\baselineskip}
\noindent
{\bf b)}  For    $GL_h(2)$  the $R$-matrix in
(\ref{1.1}) is
the solution of the Yang-Baxter equation given by
\begin{equation}\label{1.5}
R_h= \left[
\begin{array}{cccc}
1 & -h & h & h^2 \\
0 & 1 & 0 & -h \\
0 & 0 & 1 & h \\
0 & 0 & 0 & 1
\end{array} \right] \;,
\; \hat{R}_h \equiv {\cal P}R_h= \left[
\begin{array}{cccc}
1 & -h & h & h^2 \\
0 & 0 & 1 & h \\
0 & 1 & 0 & -h \\
0 & 0 & 0 & 1
\end{array} \right] \;,   \; {\cal P}R_h{\cal P}=R_h^{-1} \;,
\end{equation}
\noindent
(or $R_{h\,12}R_{h\,21}=I$, triangularity condition)
for which (\ref{1.1}) gives
\begin{equation}\label{1.6}
\begin{array}{lll}
{} [a,b]= h(\xi-a^2) \quad, & \quad [a,c]=hc^2 \quad,&\quad
[a,d]= hc(d-a) \quad, \\
{}   [b,c]=h(ac+cd) \quad,  & \quad
 [b,d]=h(d^2-\xi) \quad, & \quad [c,d]=-hc^2 \quad
\end{array}
\end{equation}
\noindent
(so that  $[a-d,c]=0$ follows),
where $\xi$ is the quadratic central element
\begin{equation}\label{1.7}
\xi \equiv det_hM = ad - cb - h cd  \quad;
\end{equation}
\noindent
setting $\xi =1$  reduces  $GL_h(2)$ to $SL_h(2)$.
The matrix $\hat{R}_h$   satisfies
\begin{equation}\label{1.7.1}
\hat{R}_h^2=I \quad, \quad (I-\hat{R}_h)(I+\hat{R}_h)=0 \quad.
\end{equation}
\noindent
It has two eigenvalues ($1$ and $-1$) and a spectral decomposition in terms
of a rank three projector $P_{h+}$ and a rank one projector $P_{h-}$
\begin{equation}\label{1.7.2}
\hat{R}_h  =P_{h+} - P_{h-} \quad, \qquad  P_{h \pm} \hat{R}_h = \pm P_{h \pm}
\quad,
\end{equation}
\begin{equation}\label{1.7.3}
P_{h+}= \frac{1}{2}(I+\hat{R}_h) \quad ,\quad P_{h-}= \frac{1}{2} (I-\hat{R}_h)
= \frac{1}{2} \left[
\begin{array}{cccc}
0 & h & -h & -h^2 \\
0 & 1 & -1 & -h \\
0 & -1 & 1 & h \\
0 & 0 & 0 & 0
\end{array} \right]  \quad.
\end{equation}
\noindent
For  $SL_h(2)$, the formulae equivalent to those in (\ref{1.4.4}) are
\begin{equation}\label{1.7.4}
\epsilon_hM^t\epsilon_h^{-1}=M^{-1} \;, \; \epsilon_h= \left(
\begin{array}{cc}
      h & 1 \\
      -1 & 0
      \end{array} \right) \;, \;  \epsilon_h^{-1}= \left(
\begin{array}{cc}
      0 & -1 \\
      1 & h
      \end{array} \right) \;, \;
P_{h- \;ij,kl}=\frac{-1}{2}\epsilon_{h\;ij}\epsilon^{-1}_{h\;kl} \;.
\end{equation}
\noindent
Using  $P_{h-}$, the deformed  determinant and its inverse, $det_hM$
and $det_hM^{-1}$,  (\ref{1.7})
are  also given by  eqs. (\ref{adet}).

\vspace{1\baselineskip}

The  quantum planes \cite{MANIN}
associated  with $SL_q(2)$ and $SL_h(2)$
are the associative algebras generated by
two elements $(x,y)\equiv X$,  the commutation properties of which
(explicitly and in $R$-matrix form) are

\vspace{0.5\baselineskip}

\noindent
{\bf a)} for $SL_q(2)$ \cite{MANIN}
\begin{equation}\label{1.8}
xy=qyx \quad  \longleftrightarrow   \quad R_qX_1X_2=qX_2X_1 \quad,
\end{equation}

\noindent
{\bf b)} for $SL_h(2)$ \cite{EOW,KAR}
\begin{equation}\label{1.9}
xy=yx+hy^2 \quad \longleftrightarrow \quad R_hX_1X_2=X_2X_1 \quad.
\end{equation}
\noindent
These  commutation
relations  are preserved under   transformations by
the corresponding quantum groups matrices\footnote{
The  $GL_q(2)$ and $GL_h(2)$ matrices also preserve, respectively, the
`$q$-symplectic' and `$h$-symplectic' metrics $\epsilon_q$
({\it or} ${\epsilon_q}^{-1}$) and ${\epsilon_h}^{-1}$.}
$M$, $X'=MX$. This invariance statement,   suitably extended  to apply
to the case of deformed Minkowski
spaces, provides the essential
ingredient for their classification.

{}From now on we shall often write $R_Q$, $P_Q$ ($Q=q,h$)
to treat both deformations simultaneously. For instance,
(\ref{1.8}) and (\ref{1.9}) may be jointly written as $R_QX_1X_2= \rho X_2X_1$,
where $\rho=(q,1)$  is the appropriate eigenvalue of $R_Q$.

\section{Deformed  Lorentz groups and associated $\qquad$ Minkowski algebras}

\indent

 As  is well known, the vector representation
$D^{\frac{1}{2}, \frac{1}{2}}= D^{\frac{1}{2},0}
\otimes D^{0, \frac{1}{2}}$  of the restricted Lorentz group
may be given by the transformation
$K'=AKA^{\dagger}$, $A \in SL(2,C)$.
The   spacetime coordinates are contained in
$K=K^{\dagger}= \sigma^{\mu} x_{\mu}$,
where $\sigma^0=I$ and $\sigma^i$ are the Pauli matrices; the time coordinate
may be identified as $x^0=\frac{1}{2}tr(K)$. Since
$detK=(x_0)^2-x^ix_i=detK'$, the correspondence
$\pm A \mapsto \Lambda \in SO(1,3)$, where
$x'^{\mu}=\Lambda^{\mu}_{\; \nu}x^{\nu}$, realizes the covering homomorphism
$SL(2,C)/Z_2=SO(1,3)$.
A first step to obtain a deformation of  the Lorentz group is to replace
the $SL(2,C)$ matrices $A$ above  by the generator
matrix $M$ of $SL_q(2)$  \cite{POWO,WA-ZPC48,SWZ-ZPC52,OSWZ-CMP}.

In general, the full determination of a deformed Lorentz group requires
the characterization of all possible commutation relations
among the generators
($a,b,c,d$) of $M$ and ($a^*,b^*,c^*,d^*$) of $M^{\dagger}$,
$M$ being a deformation of $SL(2,C)$.
The $R$-matrix form of these may be   expressed
in full generality by
\begin{equation}\label{2.3}
\begin{array}{ll}
R^{(1)}M_1M_2=M_2M_1R^{(1)} \quad, & \quad
M_1^{\dagger}R^{(2)}M_2=M_2R^{(2)}M_1^{\dagger} \quad, \\
M_2^{\dagger}R^{(3)}M_1=M_1R^{(3)}M_2^{\dagger} \quad, & \quad
R^{(4)}M_1^{\dagger}M_2^{\dagger}=M_2^{\dagger}M_1^{\dagger}R^{(4)} \quad,
\end{array}
\end{equation}
\noindent
where $R^{(3) \,\dagger}=R^{(2)}={\cal P}R^{(3)}{\cal P}$
(or `reality' condition\footnote{This reality condition
can be given in a more general
form
$R^{(3) \,\dagger}= \tau {\cal P}R^{(3)}{\cal P}$ for
$\vert \tau  \vert =1$, however this phase factor can be eliminated by
the redefinition $R^{(3)} \rightarrow \tau^{1/2} R^{(3)}$ (cf. \cite{WOZA}).}
for  $R^{(3)}$)
and   $R^{(4)}=R^{(1) \,\dagger}$ or
$R^{(4)}=({\cal P}R^{(1) \, -1}{\cal P})^{\dagger}$
since the first eq. in (\ref{2.3}) is invariant
under the exchange  $R^{(1)} \leftrightarrow
{\cal P}R^{(1) \, -1}{\cal P}$.

Eqs. (\ref{2.3}), which also follow (see {\it e.g.} \cite{PRAGA})
from the bi-spinor (dotted and undotted) description of `quantum' spacetime
in terms of a deformed
 $K$,  will be taken as the starting
point for the classification of the deformed Lorentz groups. In it, the
matrix $R^{(1)}$ characterizes  the appropriate
 deformation of the  $SL(2,C)$
group ($R^{(1)}=R_Q$), $R^{(2)}$ (or $R^{(3)}$) defines how the elements of $M$
and
$M^{\dagger}$
commute and it is not {\it a priori} fixed (but it must satisfy consistency
relations with $R^{(1)}$, see eq. (\ref{2.4})
below) and $R^{(4)}$ gives the commutation
relations for the complex conjugated generators contained in
$M^{\dagger}$.
The  specification of
the deformed  Lorentz group will be completed  by  the commutation
properties of the generators with their complex conjugated ones {\it i.e.},
by
the determination
of $R^{(2)}=R^{(3)\; \dagger}$.

The commutation relations of the deformed Lorentz group algebra generators
(entries of $M$ and $M^{\dagger}$) are given by eqs. (\ref{2.3}).
The consistency of these relations is assured if
$R^{(1)}$ (and $R^{(4)}$)
obey the Yang-Baxter equation (YBE) and $R^{(3)}$ and $R^{(2)}$ satisfy
the mixed consistency equations \cite{FTUV94-21,FM}
\begin{equation}\label{2.4}
R^{(1)}_{12}R^{(3)}_{13}R^{(3)}_{23}=R^{(3)}_{23}R^{(3)}_{13}R^{(1)}_{12}
\quad, \quad
R^{(4)}_{12}R^{(2)}_{13}R^{(2)}_{23}=R^{(2)}_{23}R^{(2)}_{13}R^{(4)}_{12}
\quad,
\end{equation}
\noindent
(these two equations are actually  the same since either
$R^{(4)}=R^{(1)\,\dagger}$
or $R^{(4)}=({\cal P} R^{(1)\,-1} {\cal P})^{\dagger}$ and
$R^{(2)}=R^{(3)\,\dagger}$). It will be convenient to notice that the
first equation, considered as an `RTT' equation, indicates that
$R^{(3)}$ is a representation of the deformed
$GL(2,C)$ group, {\it i.e.},
the matrix $R^{(3)}$ provides a 2$\times$2 representation of
the entries $M_{ij}$
of the generator matrix $M$: $(M_{ij})_{\alpha \beta}=
R^{(3)}_{i \alpha , j \beta}$. Thus, $R^{(3)}$ may be seen as a
matrix in which the 2$\times$2 blocks satisfy
among themselves the same commutation relations that the entries of $M$,
$$
R^{(3)}= \left[ \begin{array}{cc}
                   A & B \\
                   C & D
                 \end{array} \right] \; \sim \;
M = \left[ \begin{array}{cc}
                   a & b \\
                   c & d
                 \end{array} \right] \;,
$$
\noindent
and  the problem of finding all possible Lorentz deformations is
equivalent to finding all possible $R^{(3)}$ matrices with 2$\times$2 block
entries satisfying (\ref{1.3}) or (\ref{1.6}) such that
${\cal P}R^{(3)}{\cal P} = R^{(3) \dagger}$ ($ \hat{R}^{(3)}
=\hat{R}^{(3) \dagger}$).

To introduce the deformed Minkowski {\it algebra} ${\cal M}^{(j)}$
associated with a deformed Lorentz group  $L^{(j)}$ (where the index $j$
refers to the different cases)  it is natural to extend
$K'=AKA^{\dagger}$ above
to the deformed case  by stating that in it the corresponding
$K$ generates  a comodule algebra for the coaction $\phi$ defined by
\begin{equation}\label{2.5}
\phi : K \longmapsto
K' = M KM^{\dagger} \;, \; K'_{is} = M_{ij} M^{\dagger}_{ls}
K_{jl} \;, \; K=K^{\dagger} \;,\; \Lambda = M \otimes M^* \;,
\end{equation}

\noindent
where it is assumed that the matrix elements of $K$,
which now do not commute among themselves, commute with those of $M$
and $M^{\dagger}$.
As in  (\ref{1.8}), (\ref{1.9}) for $q$-two-vectors (rather, two-{\it spinors})
we now demand that the commuting
properties of the entries  of $K$  are preserved by
(\ref{2.5}).  The use of covariance arguments
to characterize the algebra generated by the entries of $K$
has been extensively used, and the resulting equations are associated
with the name of reflection equations \cite{K-SKL,KS} or, in a
more general setting,  braided algebras
\cite{MAJ-LNM,SMBR2}  of which the former constitute the `algebraic sector'
(for an introduction to braided geometry see
\cite{SMB});
similar equations were also early introduced  in \cite{FM}.
Let us now  extend the arguments given in
\cite{FTUV94-21} to classify  the
deformed Lorentz  groups and their associated Minkowski algebras in an
unified way.

This is achieved by describing the
commutation properties of the entries of the hermitian  matrix $K$ generating
a possible Minkowski algebra ${\cal M}$ by means of a
general reflection equation of the form
\begin{equation}\label{2.6}
R^{(1)} K_{1} R^{(2)} K_{2} = K_{2} R^{(3)} K_{1} R^{(4)}\quad,
\end{equation}

\noindent
where the $R^{(i)}$ matrices ($i=1,...,4$) are those
introduced in (\ref{2.3}).
Indeed, writing equation (\ref{2.6}) for $K'=MKM^{\dagger}$,
it follows that the invariance of the commutation properties of
$K$  under the associated deformed Lorentz transformation
(\ref{2.5}) is achieved if relations (\ref{2.3}) are satisfied.

The deformed  Minkowski length and metric, invariant under
a Lorentz transformation
(\ref{2.5}) of $L^{(j)}$, is defined through the quantum determinant of $K$.
Since the two matrices $\hat{R}^{(1)}={\cal P}R_Q$  have  spectral
decompositions ((\ref{1.4.2}), (\ref{1.7.2}))
with a rank three projector $P_{Q+}$ and a rank one projector $P_{Q-}$, and
the determinants of $M$, $M^{\dagger}$ are central
(eqs. (\ref{adet}), (\ref{1.7})), the $Q$-{\it deformed} and {\it invariant}
(under (\ref{2.5}))
determinant of the 2$\times$2 matrix $K$ may now be
given by
\begin{equation}\label{gdet}
(det_{Q}K)P_{Q-}P_{Q-}^{\dagger}
= - {\rho} P_{Q-} K_1 \hat{R}^{(3)}K_1 P_{Q-}^{\dagger} \quad.
\end{equation}
\noindent
It is easy to check that $(P_{Q-}P_{Q-}^{\dagger})^2=
\left( \frac{\omega_Q}{\vert [2]_{\rho} \vert} \right)^2
P_{Q-}P_{Q-}^{\dagger}$, where
$\omega_q= \vert q \vert + \vert q^{-1} \vert$, $\omega_h=2+h^2$
and $[2]_1=2$. In eq. (\ref{gdet}),
the subindex $Q$ in $det_QK$ indicates that
it depends on
$q$ or $h$ (or on other parameters on which  $R^{(3)}$ may depend)
and  ${\rho}$ ($=(q,1)$ as before) has been added by convenience.
Since $\hat{R}^{(3)}$ and $K$ are
hermitian, $det_QK$ is real (if $\rho$ is not real it may be factored out).
We stress that the above formula provides a general  expression for a central
(see below) quadratic
element which  constitutes  the {\it deformed Minkowski length}
for {\it all}  deformed spacetimes ${\cal M}^{(j)}$.

Similarly, it is possible to write in general the invariant scalar product
of {\it contravariant} (transforming  as the matrix $K$, eq.  (\ref{2.5}))
and {\it covariant} (transforming by
$Y \mapsto Y'=(M^{\dagger})^{-1}YM^{\dagger}$)
matrices (four-vectors)  as the quantum trace of a matrix product
\cite{FTUV94-21} (cf. \cite{FRT1}).  In the present general case, the
deformed trace of a matrix $B$ is defined by
\begin{equation}\label{gtr}
tr_{Q}(B):=tr({\cal D}_QB) \quad , \quad
{\cal D}_Q={\rho}^2tr_{(2)}({\cal P}(( (R_Q)^{t_1})^{-1})^{t_1}) \quad,
\end{equation}
\noindent
where $tr_{(2)}$ means trace in the second space.
This deformed trace is invariant under the quantum group coaction
$B \mapsto MBM^{-1}$ since the expression of ${\cal D}_Q$ above guarantees
that ${\cal D}_Q^t=M^t{\cal D}_Q^t(M^{-1})^t$ is fulfilled.
In particular, the ${\cal D}_Q$ matrices  for $R_q$ and $R_h$ are found
to be
\begin{equation}\label{D}
{\cal D}_q= \left( \begin{array}{cc}
        q^{-1} & 0 \\
          0   &  q
\end{array} \right) \quad, \quad {\cal D}_h=\left( \begin{array}{cc}
        1 & -2h \\
          0   &  1
\end{array} \right) \quad.
\end{equation}
\noindent
Let us now find the expression of the metric tensor. Using
$\epsilon_Q$ (cf. (\ref{1.4.4}), (\ref{1.7.4}))
$(P_{Q-})_{ij,kl}=- \frac{1}{[2]_{\rho}}\epsilon_{Q\;ij}\epsilon^{-1}_{Q\;kl}$
and  ${\cal D}_Q=
-\epsilon_Q(\epsilon_Q^{-1})^t$
(${\cal D}_Q^t=M^t{\cal D}_Q^t(M^{-1})^t$ now follows from
$\epsilon_QM^t\epsilon_Q^{-1}=M^{-1}$,
eqs. (\ref{1.4.4}) and (\ref{1.7.4})). The covariant $K^{\epsilon}_{ij}$
vector is
\begin{equation}\label{Reps}
K^{\epsilon}_{ij}=\hat{R}^{\epsilon}_{Q \, ij,kl}K_{kl} \quad, \quad
\hat{R}^{\epsilon}_Q  \equiv
(1 \otimes (\epsilon_Q^{-1})^t)
\hat{R}^{(3)} (1 \otimes (\epsilon_Q^{-1})^{\dagger}) \quad ,
\end{equation}
\noindent
from which  follows that the general Minkowski length and
metric is given by
\begin{equation}
l_Q \equiv det_QK = \frac{\rho}{\omega_Q}tr_QKK^{\epsilon} \equiv
\rho^2 g_{Q\, ij,kl} K_{ij}K_{kl} \quad , \quad
g_{Q\, ij,kl}= \frac{\rho^{-1}}{\omega_Q}
{\cal D}_{Q\, si}\hat{R}^{\epsilon}_{Q \, js,kl}\;.
\end{equation}
\noindent
This concludes the unified description of all cases. Let us now look at their
classification and specific properties.

\section{Characterization of the  Lorentz deformations}

\indent

First we use the reality
condition $R^{(3)\, \dagger}=
{\cal P}R^{(3)}{\cal P}$ to reduce the number of independent
parameters in $R^{(3)}$. It implies
\begin{equation}\label{4.1}
R^{(3)} \equiv \left[ \begin{array}{cc}
                 A&B \\
                 C&D
                \end{array}  \right] \equiv  \left[
\begin{array}{cccc}
a_{11} & a_{12} & b_{11} & b_{12} \\
a_{21} & a_{22} & b_{21} & b_{22}  \\
c_{11} & c_{12} & d_{11} & d_{12}  \\
c_{21} & c_{22} & d_{21} & d_{22}
\end{array} \right] =  \left[
\begin{array}{cccc}
a_{11} & a_{12} & a_{21}^* & b_{12} \\
a_{21} & a_{22} & b_{21} & b_{22}  \\
a_{12}^* & c_{12} & a_{22}^* & c_{22}^* \\
b_{12}^* & c_{22} & b_{22}^* & d_{22}
\end{array} \right] \;,
\end{equation}
\noindent
where $a_{11}$, $d_{22}$, $b_{21}$, $c_{12}$ are real numbers and the
rest are complex.

\vspace{0.5\baselineskip}

\noindent
{\bf a) Deformed Lorentz groups associated with $SL_q(2)$ }

\vspace{0.5\baselineskip}

Let now $M \in SL_q(2)$ and $R^{(1)}=R_q$, eq. (\ref{1.2}).
The problem of finding the $q$-Lorentz groups associated with the
standard deformation is now  reduced to obtaining all
matrices $R^{(3)}$ satisfying (\ref{2.4}). This means that
the 2$\times$2 matrices  $A,B,C,D$ in (\ref{4.1})  must  satisfy the
commutation relations in  (\ref{1.3}). This implies that
(see \cite{WOZA})   $B^2=C^2=0$,
$AD \sim I_2$ and that either $B$ or $C$ are zero. Now

\vspace{0.5\baselineskip}
\noindent
a1) $B=0$ gives  $\; R^{(3)}=\left[
\begin{array}{cccc}
a_{11} & 0 & 0 & 0  \\
0  & a_{22} & 0  & 0   \\
0  & c_{12} & a_{22}^* & 0 \\
0  & 0 & 0 & d_{22}
\end{array} \right]$ $\;$ with  $a_{11}$, $d_{22}$, $c_{12}$ $\in R\;$ .

\vspace{0.5\baselineskip}
\noindent
{}From  $AD \sim I_2$ it is easy to see  (fixing first $a_{11}=1$) that
$d_{22}=a_{22}^*/a_{22}$;  its reality then implies  $d_{22}= \pm 1$,
$d_{22}= 1$ when $a_{22} \in R$ and  $d_{22}= -1$ for   $a_{22} \in iR$.
The relation $AC=qCA$ forces   $a_{22}=q^{-1}$ or $c_{12}=0$.

\vspace{0.5\baselineskip}
\noindent
a2) $C=0$ gives  $\; R^{(3)}=\left[
\begin{array}{cccc}
a_{11} & 0 & 0 & 0  \\
0  & a_{22} & b_{21}  & 0   \\
0  & 0 & a_{22}^* & 0 \\
0  & 0 & 0 & d_{22}
\end{array} \right]$ $\,$ with  $a_{11}$, $d_{22}$,
$b_{21}$ $\in R$, $a_{11}=1$;

\vspace{1\baselineskip}
\noindent
as in the previous case, $d_{22}= \pm 1$ and  $a_{22} \in R$ for
$d_{22}= 1$    and $a_{22} \in iR$ for $d_{22}= -1$. Analogously, from
$AB=qBA$ one obtains that  $b_{21}=0$ or $a_{22}=q$.

\vspace{1\baselineskip}

Thus,  the solutions for $R^{(3)}$ are the following
\begin{eqnarray}
 & & R^{(3)}=\left[
\begin{array}{cccc}
q & 0 & 0 & 0  \\
0  & 1 & 0  & 0   \\
0  & r & 1 & 0 \\
0  & 0 & 0 & q
\end{array} \right] \quad, \quad   \begin{array}{l}
                              q \in R  \; ,\\
                              r \in R \;,
                                   \end{array} \,   \label{4.3}\\
 & & R^{(3)}=\left[
\begin{array}{cccc}
1 & 0 & 0 & 0  \\
0  & t & 0  & 0   \\
0  & 0 & \pm t & 0 \\
0  & 0 & 0 & \pm 1
\end{array} \right] \quad, \quad   \begin{array}{l}
                              + \; \mbox{for} \; t \in R  \; ,\\
                              - \; \mbox{for} \; t \in iR \;,
                                   \end{array} \,  \label{4.2}\\
 & & R^{(3)}=\left[
\begin{array}{cccc}
q^{-1} & 0 & 0 & 0  \\
0  & 1 & r  & 0   \\
0  & 0 & 1 & 0 \\
0  & 0 & 0 & q^{-1}
\end{array} \right] \quad, \quad   \begin{array}{l}
                              q \in R  \; ,\\
                              r \in R \;,
                                   \end{array} \,  \label{4.4}\\
 & & R^{(3)}=\left[
\begin{array}{cccc}
1 & 0 & 0 & 0  \\
0  & q^{-1} & 0  & 0   \\
0  & r & -q^{-1} & 0 \\
0  & 0 & 0 & -1
\end{array} \right] \quad, \quad   \begin{array}{l}
                              q \in iR  \; ,\\
                              r \in R \;,
                                   \end{array} \,  \label{4.5}\\
 &  & R^{(3)}=\left[
\begin{array}{cccc}
1 & 0 & 0 & 0  \\
0  & q & r  & 0   \\
0  & 0 & -q & 0 \\
0  & 0 & 0 & -1
\end{array} \right] \quad, \quad   \begin{array}{l}
                              q \in iR \; ,\\
                              r \in R \;,
                                   \end{array} \, \label{4.6}
\end{eqnarray}
\noindent
{\it Remarks:}

\noindent
- Notice that, as anticipated, the $Q$-`determinant' of all these
$R^{(3)}$ matrices, computed as $det_Q M$, is a scalar -hence commuting-
2$\times$2 matrix.

\noindent
- $R_q^{\dagger}={\cal P}R_q{\cal P}$ iff  $q \in R$.
Hence, $R^{(4)}_{12}=R_{q21}$ or $R^{-1}_{q12}$. Thus
$\tilde{M}\equiv  (M^{-1})^{\dagger}$  provides
a second copy of $SL_q(2)$, since then
$R_q\tilde{M}_1\tilde{M}_2=\tilde{M}_2\tilde{M}_1R_q$.

\noindent
- The case (\ref{4.3})  for $r=q-q^{-1}=\lambda$ ($R^{(3)}=R_q$)
is the quantum Lorentz group of  \cite{WA-ZPC48,SWZ-ZPC52}
($L_q^{(1)}$  in the notation of \cite{FTUV94-21}).
If $r \neq \lambda$ we obtain a `gauged'
version of it:
$R^{(3)}=e^{ \alpha \sigma^3_2} R_q e^{ - \alpha \sigma^3_2}$
($r=\lambda e^{2 \alpha}$), where the subindex in $\sigma^3_2$ refers to
the second space.

\noindent
- The matrix  (\ref{4.2})   for $t=1$ and $q \in R$ corresponds to
$L_q^{(2)}$  in \cite{FTUV94-21}.

\noindent
- The calculations leading to (\ref{4.2})-(\ref{4.6})
require assuming  $q^2 \neq 1$.
However,  the  solutions for $q \in R$ are also valid in the limit $q=1$
(see \cite{WOZA}); in this limit ($R^{(1)}=R^{(4)}=I_4$),
the case (\ref{4.2})  gives
the deformed Lorentz group (twisted) of \cite{CHA-DEM}.
For $q=-1$, additional solutions appear  and, although
we shall not discuss these particular cases (see \cite{WOZA}), the
associated Minkowski
algebras may be obtained as in the general $q$ case.

\noindent
- These results coincide with the classification in \cite{WOZA}:
the solutions (\ref{4.2}) correspond to eqs. (13) and (14) in \cite{WOZA};
similarly, (\ref{4.3}), (\ref{4.4}), (\ref{4.5}) and (\ref{4.6})
correspond to (74) ($q$ real), (15), (74) ($q$ imaginary) and (16)  in that
reference.

\vspace{0.5\baselineskip}

\noindent
{\bf b) Deformed Lorentz groups associated with $SL_h(2)$ }

\vspace{0.5\baselineskip}

Let now $R^{(1)}=R_h$, eq. (\ref{1.5}).  For $h$
imaginary, $h \in iR$, the matrix $R_h$ satisfies the reality condition
$R^*_h=R_h^{-1}$ ($={\cal P}R_h {\cal P}$); this means that
$\tilde{M} \equiv M^*$ defines a second copy of $SL_h(2)$ since
$R_hM^*_1M^*_2=M^*_2M^*_1R_h$. The value of $h \in C \backslash \{ 0 \}$,
however, is not important. Indeed, quantum groups related with two
different values of $h \in C$
are equivalent and  their $R$ matrices  are related by  a similarity
transformation\footnote{Quantum groups associated with $R_h$ and $R_{h=1}$
are related by a similarity transformation defined by the 2$\times$2
matrix $S= diag(h^{-1/2} , h^{1/2})$:  $R_{h=1}=
(S \otimes S) R_h (S \otimes S)^{-1}$.}; thus, we can take
$h \in R$ or even $h=1$.

Since the entries of $M$ satisfy (\ref{1.6}), the  2$\times$2  blocks
in  $R^{(3)}$  (eq. (\ref{4.1})) will satisfy now these  commutation relations.
This leads to (see \cite{WOZA})  $C=0$
so that, taking the $h$-`determinant' of $R^{(3)}$ equal $I_2$,
the set of commutation
relations reduces to
\begin{equation}\label{4.7}
AD=I_2 \quad, \qquad [A,B]=h(I_2 -A^2) \quad.
\end{equation}
\noindent
Using them in (\ref{4.1}) the following solutions for $R^{(3)}$
are found ($h \in R$)
\begin{eqnarray}
& & R^{(3)}=\left[
\begin{array}{cccc}
1 & 0 & 0 & 0  \\
0  & 1 & r  & 0   \\
0  & 0 & 1 & 0 \\
0  & 0 & 0 & 1
\end{array} \right] \quad, \quad   r \in R \quad ,   \, \label{4.8}\\
 & & R^{(3)}=\left[
\begin{array}{cccc}
1 & 0 & -h & 0  \\
-h  & 1 & r  & h   \\
0  & 0 & 1 & 0 \\
0  & 0 & h & 1
\end{array} \right] \quad, \quad   \begin{array}{l}
                              h \in R  \; ,\\
                              r \in R \;.
                                   \end{array} \, \label{4.9}
\end{eqnarray}
\noindent
{\it Remarks:}

\noindent
-  In (\ref{4.9}), for $r=h^2$ we have
$R^{(3)}=({\cal P}R_h{\cal P})^{t_2}$. However,
the parameter $r$ can be removed with an appropriate change of basis
provided $h \neq 0$.  For $h=0$, this is not possible and
constitutes a different case, eq. (\ref{4.8}). This case is another example
where the non-commutativity    is solely due to $R^{(3)} \neq I_4$.

\noindent
- The cases (\ref{4.8}), (\ref{4.9}) correspond to (20)
and (21) [cf. (78)]  in \cite{WOZA}.

\section{Minkowski algebras: classification and properties }

\indent

We now present here, in  explicit form, the commutation relations
for the generators of the deformed Minkowski spacetimes; they follow
easily
from  (\ref{2.6}) using the previous $R^{(3)}$ matrices.
We saw in (\ref{2.3}) that $R^{(3) \,\dagger}=R^{(2)}={\cal P}R^{(3)}{\cal P}$
and   $R^{(4)}=R^{(1) \,\dagger}$ or
$R^{(4)}=({\cal P}R^{(1) \, -1}{\cal P})^{\dagger}$
(these two possibilities are the same for $Q=h$). Clearly, eq. (\ref{2.6})
allows for a factor in one side without impairing its invariance
properties. This factor may be selected with the (natural) condition
that the resulting Minkowski algebra does not contain generators
$\alpha, \beta,...$, with the Grassmann-like property
$\alpha^2=\beta^2= ...=0$. In terms of $P_{Q+}$, this tantamount to
requiring that $P_{Q+}K_1 \hat{R}^{(3)}K_1P_{Q+}^{\dagger}$ must
be non-zero.
This leads to (cf. (\ref{2.6}))
the equations
\begin{equation}\label{RE2}
R_QK_1R^{(2)}K_2= \pm K_2R^{(3)}K_1R_Q^{\dagger} \quad
(+ \; \mbox{for} \; q,h \in R \; , \;
- \; \mbox{for} \; q\in iR )\;.
\end{equation}
\noindent
In the $q$-case we might also consider
$R^{(4)}=({\cal P}R^{(1) \, -1}{\cal P})^{\dagger}$. However
using  Hecke's  condition for   $R^{(1)}$  it is seen that this leads
to the same algebra as (\ref{RE2}) with the restriction $det_qK=0$,
so that this case may be considered as included in the previous one.

An important ingredient  is the centrality of the $Q$-determinant
(\ref{gdet}), $(det_QK)K=K(det_QK)$, since it will correspond to the Minkowski
length. Using twice (\ref{RE2}) we find the following commutation property
for three $K$ matrices
\begin{equation}
R_{Q\,13}R_{Q\,23}K_1R^{(2)}_{12} K_2 R^{(2)}_{13}R^{(2)}_{23}K_3=
K_3 R^{(3)}_{13}R^{(3)}_{23}K_1R^{(2)}_{12}K_2 R^{\dagger}_{Q\,13}
R^{\dagger}_{Q\,23}\;.
\end{equation}
\noindent
Multiplying from the right by ${\cal P}_{12}P_{Q-\;12}^{\dagger}$
and by $P_{Q-\;12}$  from the left and using that $R_Q$ and
$R^{(3)}$ represent $GL_Q(2)$  and hence have a central
$Q$-`determinant'
represented  by a scalar 2$\times$2 matrix we get
\begin{equation}\label{centr}
(det_QR_Q)(det_QR^{(3)})^{\dagger} \,(det_QK)K=
(det_QR_Q)^{\dagger}(det_QR^{(3)}) \, K(det_QK) \quad,
\end{equation}
\noindent
The scalar  $det_QR^{(i)}$ matrices always cancel out in the cases below
($det_qR_q=qI_2$ and $det_hR_h=I_2$)
assuring the centrality of $det_QK$ (as it may
be checked by direct  computation).

\vspace{0.5\baselineskip}

\noindent
{\bf a)  $q$-Minkowski spaces associated with $SL_q(2)$}

\vspace{0.5\baselineskip}

\noindent
1) Let us consider the case (\ref{4.3}) for
$r= \lambda$ ({\it i.e.}, $R^{(3)}=R_q$, $q$ real). The commutation relations
for the entries of $K$={\footnotesize $   \left( \begin{array}{cc}
                        \alpha & \beta \\
                        \gamma & \delta
                        \end{array}  \right)$}
are
\begin{equation}\label{5.4}
\begin{array}{lll}
\alpha \beta = q^{-2} \beta \alpha \quad , &
[\delta, \beta ] = q^{-1} \lambda \alpha \beta \quad ,&

\alpha \gamma = q^{2} \gamma \alpha \quad, \\

{} [\beta , \gamma] = q^{-1} \lambda (\delta - \alpha ) \alpha \quad ,&
[\alpha , \delta ] = 0 \quad ,&
 [\gamma , \delta ]
= q^{-1} \lambda \gamma \alpha \quad ;
\end{array}
\end{equation}
\noindent
they  characterize the algebra ${\cal M}_{q}^{(1)}$
(\cite{WA-ZPC48}-\cite{OSWZ-CMP}; see also
\cite{SMBR2,AKR,MeyerM,FTUV94-21}).
The Minkowski length is given by  (\ref{gdet}),
\begin{equation}\label{5.5}
det_qK=  \alpha \delta - q^2 \gamma \beta \quad.
\end{equation}
\noindent
If $r \neq \lambda$, the commutation relations are slightly different;
this, however, corresponds only to an appropriate election of the basis
(`gauged' version of this Minkowski space).

\vspace{1\baselineskip}
\noindent
2) Let  $R^{(3)}$ be given by eq. (\ref{4.2}). The centrality of the
$q$-determinant implies that $q$ and $t$ are both real or both imaginary.
The commutation relation for the entries of $K$
and the $q$-Minkowski length (eq. (\ref{gdet}))  are (the
sign $+$ is for $q,t \in R$ and the $-$ for $q,t \in iR$)
\begin{equation}\label{5.1}
\begin{array}{lll}
q \alpha \beta = \pm t \beta \alpha   \quad, &
\; t \alpha \gamma = \pm q \gamma \alpha \quad, & \;
\alpha \delta =  \delta \alpha \quad, \\
{}[\beta, \gamma]= \pm t \lambda \alpha \delta \quad, & \;
\beta \delta = \pm qt \delta \beta \quad, &
\; \delta \gamma = \pm qt \gamma \delta \quad;
\end{array}
\end{equation}
\begin{equation}\label{5.3}
det_{q,t}K=   \frac{q+ q^{-1}}{q \pm q^{-1}}
(- q \gamma \beta   \pm t \alpha \delta)  \quad.
\end{equation}
\noindent
{\it Remarks:}

\noindent
- For $t=1$, these commutation relations correspond to the Minkowski
algebra ${\cal M}^{(2)}_q$ \cite{WA-ZPC48,MajEucl,FTUV94-21} which
is isomorphic to the quantum algebra\footnote{The Minkowski space of
\cite{VD} is also a $GL_q(2)$-like space, but different from the above.}
$GL_q(2)$.

\noindent
- For $q=1$ and $t$ real, we get  the Minkowski space
obtained in \cite{CHA-DEM}
(denoted ${\cal M}^{(3)}$  in \cite{FTUV94-21}).
This algebra and the corresponding deformed Poincar\'e algebra
have been shown to be \cite{LRT} a simple transformation
(twisting) of the classical one. As a result, it is possible
to remove the non-commuting character of the
entries of $K$ \cite{AKRREL}.

\vspace{1\baselineskip}

\noindent
3) Let us  take $R^{(3)}$ as in eq. (\ref{4.4}) for
$r$=$- \lambda$ ($R^{(3)}={\cal P}R^{-1}_q{\cal P}$).
Then,
\begin{equation}\label{5.6}
\begin{array}{lll}
{} [ \alpha, \beta]=q \lambda \beta \delta \quad, &
\; [ \alpha, \gamma ]= - q \lambda \delta \gamma \quad, & \;
[ \alpha, \delta ]= 0 \quad, \\
{} [ \beta, \gamma ]
=q \lambda (\alpha - \delta) \delta \quad,  & \;
\beta \delta = q^2 \delta \beta \quad, & \;
\gamma \delta = q^{-2} \delta \gamma \quad;
\end{array}
\end{equation}
\begin{equation}\label{5.6.1}
det_qK=q^2 \alpha \delta - \beta \gamma \quad.
\end{equation}
\noindent
This algebra may also be identified with the algebra of spacetime
derivatives in \cite{OSWZ-CMP} (see also \cite{AKR}).

\vspace{1\baselineskip}

\noindent
4) Let  $R^{(3)}$ be now  given by (\ref{4.5}).
The Minkowski algebra and the central length are given by
\begin{equation}\label{5.8}
\begin{array}{lll}
\alpha \beta = -q^{-2} \beta \alpha \quad ,\; &
\delta \beta + \beta \delta  = r \alpha \beta \quad ,\;&
\alpha \gamma = - q^{2} \gamma \alpha \quad, \\
{} [\beta , \gamma] = - q^{-1} \lambda  \delta  \alpha  +r \alpha^2 \quad ,\;
& [\alpha , \delta ] = 0 \quad ,\; &
\gamma \delta + \delta \gamma
= r  \gamma \alpha \quad ,
\end{array}
\end{equation}
\begin{equation}\label{5.8.1}
det_qK=  \frac{-q [2]}{\lambda}(q^{-2} \alpha \delta +  \gamma \beta )\quad.
\end{equation}
\noindent

\vspace{1\baselineskip}

\noindent
5) Finally, let $R^{(3)}$ be as  in eq. (\ref{4.6}). Then,
\begin{equation}\label{5.9}
\begin{array}{lll}
\alpha \beta + \beta \alpha = -r  \beta \delta \quad, &
\;  \alpha \gamma + \gamma \alpha = - r \delta \gamma \quad, & \;
[ \alpha, \delta ]= 0 \quad, \\
{} [ \beta, \gamma ]
 = - q \lambda \alpha \delta + r \delta^2   \quad,  & \;
\beta \delta = - q^2 \delta \beta \quad, & \;
\gamma \delta = - q^{-2} \delta \gamma \quad;
\end{array}
\end{equation}
\begin{equation}\label{5.9.1}
det_qK= \frac{-q [2]}{\lambda} ( q^2 \alpha \delta + \beta \gamma) \quad.
\end{equation}
\noindent

\vspace{0.5\baselineskip}

\noindent
{\bf b) Deformed Minkowski spaces associated with $SL_h(2)$}

\vspace{0.5\baselineskip}

\noindent
1) Let $R^{(3)}$ be given first by  eq. (\ref{4.8})
and let $R^{(1)}=R_h$, eq. (\ref{1.5}).
Using (\ref{RE2}) with the plus sign and  (\ref{gdet})
we find  ($h$ real)
\begin{equation}\label{5.7}
\begin{array}{ll}
{} [ \alpha, \beta ]= -h \beta^2 -r \beta \delta + h \delta \alpha
-h \beta \gamma + h^2 \delta \gamma \;,&
\quad [\alpha, \delta]= h( \delta \gamma - \beta \delta) \;,\\
{} [\alpha, \gamma]= h \gamma^2 + r \delta \gamma - h \alpha \delta
+h \beta \gamma - h^2 \beta \delta \;, & \quad
[\beta, \delta]=  h \delta^2 \;,\\
{} [ \beta, \gamma]= h \delta ( \gamma + \beta) + r \delta^2 \;,& \quad
[ \gamma, \delta]= - h \delta^2 \;;
\end{array}
\end{equation}
\begin{equation}\label{5.7.1}
det_{h} K= \frac{2}{h^2+2}(\alpha \delta - \beta \gamma + h \beta \delta)
\quad.
\end{equation}

\vspace{1\baselineskip}

\noindent
2) Let $R^{(3)}$ be given now by  eq. (\ref{4.9}) with $r=0$.
In this case,
\begin{equation}\label{5.10}
\begin{array}{ll}
{} [ \alpha, \beta ]= 2h \alpha \delta + h^2 \beta \delta  \quad , &
\quad [\alpha, \delta]= 2h (\delta \gamma -  \beta \delta ) \quad,\\
{} [\alpha, \gamma]= -h^2 \delta \gamma - 2h \delta \alpha \quad, & \quad
[\beta, \delta]=  2 h \delta^2 \quad , \\
{} [ \beta, \gamma]= 3h^2 \delta^2  \quad ,& \quad
[ \gamma, \delta]= - 2h \delta^2 \quad ;
\end{array}
\end{equation}
\begin{equation}\label{5.11}
det_{h} K=\frac{2}{h^2+2}(\alpha \delta - \beta \gamma + 2h \beta \delta)
\quad.
\end{equation}

\vspace{0.5\baselineskip}

\noindent
{\bf c) Final remarks}

\vspace{0.5\baselineskip}

\noindent
For all the $Q$-spacetime algebras, time may be defined as proportional to
$tr_QK$ ($=2x^0$ in the undeformed case). The time generator obtained
in this way is central only for ${\cal M}_q^{(1)}$
\cite{WA-ZPC48}-\cite{OSWZ-CMP} and for the Minkowski algebra (\ref{5.6})
(in fact, they are isomorphic: the entries of the
covariant vector $K^{\epsilon}$ for ${\cal M}_q^{(1)}$
satisfy the commutation relations  (\ref{5.6}) \cite{FTUV94-21}).

The differential calculus on all the above Minkowski spaces
may be easily discussed now along the lines of \cite{FTUV94-21,AKR};
one could also investigate the r\^ole played in it by the contraction relating
\cite{AKS} the $q$- and  $h$-deformations.
To conclude, let us mention that the additive braided group
structure \cite{MAJ-LNM}-\cite{SMB}
of all these algebras,
may be easily found. It suffices to impose that
eq. (\ref{RE2})  is also satisfied by the sum $K'+K$ of two
copies $K$ and $K'$. Using  Hecke's condition
($R_{Q12}=R_{Q21}^{-1}+(\rho - \rho^{-1}){\cal P}$) this
gives
\begin{equation}\label{braid}
R_QK_1'R^{(2)}K_2= \pm K_2R^{(3)}K_1'({\cal P}R_Q^{\dagger}{\cal P})^{-1} \quad
(+ \; \mbox{for} \; q,h \in R \; , \;
- \; \mbox{for} \; q\in iR )\;,
\end{equation}
\noindent
which is clearly preserved by (\ref{2.5}); for ${\cal M}_q^{(1)}$,
it reproduces the result of \cite{MeyerM}.

\vspace{1\baselineskip}

\noindent
{\bf Acknowledgements:} This article has been partially supported by a
research grant from the CICYT, Spain. The authors are indebted to P.P. Kulish
for many useful discussions.

\end{document}